\documentclass{amsart}

\usepackage{amssymb,latexsym,amsmath,amsfonts}
\input xy
\xyoption{all}

\newcommand{\cal}{\mathcal}

\makeatletter
\renewcommand{\subsection}{\@startsection{subsection}{2}{0mm}{-2mm}{-2mm}{\bf\normalsize}}
\def\sbsn{\subsection{\hspace{-3mm}}}

\def\sbsnt#1{\subsection{#1}}
\makeatother

\newtheorem{formula}{}[section]
\newtheorem{definition}[formula]{Definition}
\newtheorem{corollary}[formula]{Corollary}
\newtheorem{remark}[formula]{Remark}
\newtheorem{lemma}[formula]{Lemma}
\newtheorem{theorem}[formula]{Theorem}

\def\thrm{\begin{theorem}}
\def\thrml#1{\begin{theorem}\label{#1}}
\def\ethrm{\end{theorem}}
\def\rmrk{\begin{remark}}
\def\rmrkl#1{\begin{remark}\label{#1}}
\def\ermrk{\end{remark}}
\def\dfntn{\begin{definition}}
\def\dfntnl#1{\begin{definition}\label{#1}}
\def\edfntn{\end{definition}}
\def\nmrt{\begin{enumerate}}
\def\enmrt{\end{enumerate}}
\def\tm#1{\item[{\rm (#1)}]}
\def\qtn{\begin{equation}}
\def\qtnl#1{\begin{equation}\label{#1}}
\def\eqtn{\end{equation}}
\def\lmm{\begin{lemma}}
\def\lmml#1{\begin{lemma}\label{#1}}
\def\elmm{\end{lemma}}
\def\crllr{\begin{corollary}}
\def\crllrl#1{\begin{corollary}\label{#1}}
\def\ecrllr{\end{corollary}}
\def\css{\begin{cases}}
\def\ecss{\end{cases}}

\def\proof{\noindent{\bf Proof}.\ }

\def\A{{\cal A}}

\def\P{{\cal P}}

\def\T{{\cal T}}
\def\X{{\cal X}}

\def\CC{{\mathbb C}}

\def\ZZ{{\mathbb Z}}
\def\Q{{\mathbb Q}}

\def\SP{{\scriptscriptstyle P}}
\def\SPG{{\scriptscriptstyle P^G}}

\DeclareMathOperator{\aut}{Aut}
\DeclareMathOperator{\AG}{AG}
\DeclareMathOperator{\AGL}{AGL}
\DeclareMathOperator{\cla}{Cla}
\DeclareMathOperator{\cyc}{Cyc}
\DeclareMathOperator{\GCD}{GCD}
\DeclareMathOperator{\GF}{GF}
\DeclareMathOperator{\GL}{GL}
\DeclareMathOperator{\fix}{Fix}
\DeclareMathOperator{\id}{id}

\DeclareMathOperator{\Irr}{Irr}
\DeclareMathOperator{\iso}{Iso}
\DeclareMathOperator{\mat}{Mat}
\DeclareMathOperator{\orb}{Orb}
\DeclareMathOperator{\PGL}{PGL}
\DeclareMathOperator{\PGaL}{P{\rm \Gamma}L}
\DeclareMathOperator{\PSL}{PSL}

\DeclareMathOperator{\reg}{reg}
\DeclareMathOperator{\rk}{rk}

\DeclareMathOperator{\sym}{Sym}
\DeclareMathOperator{\tr}{tr}

\def\bul{\hfill\vrule height .9ex width .8ex depth -.1ex}
\def\bull{\hfill\vrule height .9ex width .8ex depth -.1ex\medskip}

\def\lg{\langle}

\def\ov{\overline}
\def\rg{\rangle}
\def\wt{\widetilde}

\def\VRT#1{*=<5mm>[o][F-]{#1}}
\def\EVR{\VRT{}}
\def\grphp#1{$\xymatrix@R=10pt@C=10pt@M=0pt@L=2pt{#1}$}

\begin{document}
\title{On Pseudocyclic Association Schemes}
\author{Mikhail Muzychuk}
\address{Netanya Academic College, Netanya, Israel}
\email{muzy@netanya.ac.il}
\author{Ilya Ponomarenko}
\address{Steklov Institute of Mathematics at St. Petersburg, Russia}
\email{inp@pdmi.ras.ru}
\date{}

\begin{abstract}
The notion of pseudocyclic association scheme is generalized to the non-commutative
case. It is proved that any pseudocyclic scheme the rank of which is much more than
the valency is the scheme of a Frobenius group and is uniquely determined up to
isomorphism by its intersection number array. An immediate corollary of this result
is that any scheme of prime degree, valency $k$ and rank at least $k^4$ is schurian.
\end{abstract}

\maketitle

\section{Introduction}
A commutative association scheme is called {\it pseudocyclic} if the multiplicity of its
non-principal irreducible character does not depend on the choice of the character~\cite{BCN}
(for a background on theory of association schemes and coherent configurations we
refer to \cite{EP09} and Section~\ref{220309a}). It can be proved that such a scheme
is {\it equivalenced}, i.e. the valencies of its non-reflexive basis relations are pairwise equal. A classical example
of a commutative pseudocyclic scheme is a cyclotomic scheme over a finite field; two other
series of such schemes were constructed in~\cite{HX06}. A motivation to study
pseudocyclic schemes is that any of them produces special $2$-designs.\medskip

In this paper we define an association scheme (not necessary commutative) to be
{\it pseudocyclic} if the ratio of the multiplicity and degree of its non-principal
irreducible character does not depend on the choice of the character. Clearly,
in the commutative case both definitions give the same concept. To formulate
one of the main results of this paper (which, in particular, shows that non-commutative
pseudocyclic schemes do exist) we make two remarks. First, as in the commutative case
we can prove (Theorem~\ref{pscyclic}) that any pseudocyclic scheme is equivalenced;
the valency of its non-reflexive basis relation is
called the {\it valency} of the scheme. Secondly, under a {\it Frobenius scheme} we
mean the scheme of a Frobenius group (in its standard permutation representation in
which the one point stabilizer coincides with the Frobenius complement). Our first
result is an immediate consequence of Theorems~\ref{190309g} and~\ref{220309h}.

\thrml{020209a}
Any Frobenius scheme is pseudocyclic. Conversely, there exists a function $f(k)$ such
that any pseudocyclic scheme of valency $k>1$ and rank at least $f(k)$ is a Frobenius
scheme.~\bull
\ethrm

The rough upper bound on the function $f(k)$ obtained in the proof is $O(k^4)$. On
the other hand, the scheme of a non-Desarguesian affine plane of order $q$ is
non-schurian (see Subsection~\ref{100209d}), and hence can not be a Frobenius
scheme. Besides, it is pseudocyclic of valency $q-1$ and rank $q+2$. Thus
$f(k)\ge k+3$.\medskip

The proof of the second part of Theorem~\ref{020209a} is based on Theorem~\ref{280109c}
giving together with Theorem~\ref{250309a} a sufficient condition for an equivalenced
scheme to be schurian.

The celebrated Hanaki-Uno theorem states that any scheme of prime degree is
pseudocyclic~\cite{HU}. Therefore as an immediate consequence of Theorem~\ref{020209a} we
have the following result.

\crllr
Any scheme of prime degree, valency $k$ and rank at least $f(k)$ is schurian.\bull
\ecrllr

One of the most important problems in association scheme theory is to determine
a scheme up to isomorphism by means of the intersection number array. For example, the
intersection number array of the scheme of a distance-regular graph is uniquely
determined by the parameters of the graph. Therefore the most part of characterizations
of the classical distance-regular graphs given in~\cite{BCN} are in fact the
characterizations of their schemes in the above sense (we refer to~\cite{EP09} for
a survey of relevant results). In this paper we prove the following theorem.

\thrml{020209b}
Any Frobenius scheme of valency $k$ and rank at least $f(k)$ is determined up to
isomorphism by its intersection number array.
\ethrm

All undefined terms and results concerning permutation groups can be found in
monographs~\cite{Pa,W64}. To make the paper as self-contained as possible we cite the
background on association schemes~\cite{Zi1} and coherent configurations~\cite{EP09}
in Section~\ref{220309a}. Section~\ref{310309a} contains the definition
of a pseudocyclic scheme, several useful results on these schemes and the proof of the first
part of Theorem~\ref{020209a} (Theorem~\ref{190309g}). In Section~\ref{040209a}
we give a brief exposition of known families of pseudocyclic schemes. In
Sections~\ref{090209g} and~\ref{310309b} under a special assumption we explicitly find
a one point extension of a pseudocyclic scheme and a one point extension of any algebraic
isomorphism from it to another scheme (Theorems~\ref{280109c} and~\ref{280109e}). Based on these results we prove
our main theorems in Section~\ref{310309c}.
Section~\ref{310309d} includes some concluding
remarks and special results concerning pseudocyclic schemes.\medskip

{\bf Notation.}
Throughout the paper $\Omega$ denotes a finite set. Set
$1_\Omega=\{(\alpha,\alpha):\ \alpha\in\Omega\}$ and $1_\alpha=1_{\{\alpha\}}$ for all
$\alpha\in\Omega$. For a relation $r\subset\Omega\times\Omega$ we set
$r^*=\{(\beta,\alpha):\ (\alpha,\beta)\in r\}$ and
$\alpha r=\{\beta\in\Omega:\ (\alpha,\beta)\in r\}$ for all $\alpha\in\Omega$.
The adjacency matrix of~$r$ is denoted by $A(r)$.
For $s\subset \Omega\times\Omega$ we set
$r\cdot s=\{(\alpha,\gamma):\ (\alpha,\beta)\in r,\ (\beta,\gamma)\in s$
for some $\beta\in\Omega\}$. If $S$ and $T$ are sets of relations, we set
$S\cdot T=\{s\cdot t:\ s\in S,\, t\in T\}$. For a permutation group $G\le\sym(\Omega)$ we
denote by $\orb(G)=\orb(G,\Omega)$ the set of $G$-orbits.

\section{Pseudocyclic schemes}\label{310309a}

A scheme $(\Omega,S)$ is called {\it pseudocyclic} if the number $m_\SP/n_\SP$ does not
depend on the choice of the central primitive idempotent $P\in\P^\#$ (see Subsection~\ref{190809a}). In the
commutative case $n_\SP=1$ for all $P\in\P$, and our definition is compatible with that
from~\cite{BCN}. Any regular scheme (not necessarily commutative) is pseudocyclic because in this case $m_\SP=n_\SP$
for all~$P$. More elaborated example of a non-commutative pseudocyclic scheme arises
from a Frobenius group with non-abelian kernel. In what follows we use the family of
such groups given in \cite[pp.187-189]{Pa}.\medskip

{\bf Example.} Let $q$ be a prime power and $n>1$ an odd integer. Set $H$ to be a subgroup
of $\GL(3,q^n)$ that consists of all matrices of the form
$$
A(a,b)=
\begin{pmatrix}
1 & a & b \\
0 & 1 & a^q\\
0 & 0 & 1\\
\end{pmatrix}
\qquad a,b\in\GF(q^n).
$$
Then the mapping $\sigma:A(a,b)\mapsto A(ca,c^{1+q}b)$ is a fixed point free automorphism
of $H$ whenever the multiplicative order of $c$ equals $(q^n-1)/(q-1)$.
So the semidirect product $G=HK$ where $K=\lg\sigma\rg$, is a Frobenius group with
non-abelian kernel $H$ and complement $K$. The natural action of~$G$ on~$H$ produces an
equivalenced scheme of degree $|H|=q^{2n}$, rank $q^{n+1}-q^n+q$ and valency
$(q^n-1)/(q-1)$. A straightforward calculation for $(q,n)=(2,3)$ shows that the adjacency
algebra of the corresponding scheme has exactly four irreducible characters the
multiplicities and degrees of which are as follows:
$$
(m_1,n_1)=(1,1),\ (m_2,n_2)=(7,1),\ (m_3,n_3)=(m_4,n_4)=(14,2).
$$
In particular, the scheme is not commutative but is pseudocyclic  because $m_\SP/n_\SP=7$
for all $P\in\P^\#$. In fact, this example is a special case of the following theorem.

\thrml{190309g}
Any Frobenius scheme is pseudocyclic. Moreover, it is commutative if and only if the
kernel of the associated Frobenius group is abelian.
\ethrm
\proof Let $(\Omega,S)$ be a Frobenius scheme, $G\leq\sym(\Omega)$ the corresponding Frobenius
group and the mapping
$$
\pi:\CC G\rightarrow\mat_\Omega(\CC)
$$
is induced by
the permutation representation of $G$ (in particular, $\pi(g)$ is the permutation matrix
of a permutation $g\in G$). It is a well-known fact that a semisimple subalgebra of
$\mat_\Omega(\CC)$ coincides with the centralizer of its centralizer (in $\mat_\Omega(\CC)$)
\cite[p.178]{CR}. {Since the algebra $\pi(\CC G)$ is semisimple and its centralizer in
$\mat_\Omega(\CC)$ equals $\CC S$, the centralizer of $\CC S$ in $\mat_\Omega(\CC)$
coincides with $\pi(\CC G)$. Therefore these two algebras have the same centre
\qtnl{280509a}
Z(\pi(\CC G))=Z(\CC S)=\CC S\cap \pi(\CC G),
\eqtn
}
and hence the same set of central primitive idempotents, say~$\P$. Moreover, one can see
that $n_{\SP}=m_{\chi_\SP}$ and $m_{\SP}=\chi_{\SP}(1)$ for all $P\in\P$ where $\chi_P$ is the
irreducible character of $G$ corresponding to the central primitive
idempotent~$P$. Since $G$ is a Frobenius group, by Theorem~\ref{230309a} this implies
that
$$
m_\SP/n_\SP=m_{\chi_\SP}/\chi_{\SP}(1)=|K|,\qquad P\in\P^\#,
$$
where $K$ is the complement of~$G$. Thus $(\Omega,S)$ is a pseudocyclic scheme.\medskip

To prove the second statement we find the dimension of the left-hand side
of~(\ref{280509a}). First, we note that the algebra $Z(\CC G)$ is spanned by the elements
$C^+=\sum_{g\in C} g$ where $C$ runs through the set of conjugacy classes of~$G$.
Moreover, if $A$ denotes the kernel of $G$, then
$$
C\cap A=\emptyset\ \Rightarrow\ \pi(C^+)\in\CC J_\Omega.
$$
(Indeed, $C$ contains $a^{-1}ca=a^{-1}a^cc$ for all $a\in A$ and $c\in C$. On
the other hand, since $c\not\in A$, the mapping $a\mapsto a^{-1}a^c$, $a\in A$,
is a bijection. Thus $CA=C$. Taking into account that $\pi(A^+)=J_{\Omega}$, we
conclude that $\pi(C^+)$ is a scalar multiple of $J_\Omega$.) This implies that
\qtnl{300509a}
\dim(\pi(Z(\CC G)))=|\cla_G(A)|
\eqtn
where $\cla_G(A)$ is the set of conjugacy classes of~$G$ contained in~$A$.
Indeed, $\pi(C^+)$ is a $\{0,1\}$ matrix. Since the sum of these matrices
with $C\in\cla_G(A)$ equals $J_\Omega$, they are linearly independent and form a
basis of $\pi(Z(\CC G))$.\medskip

Next, the group $K$ acts semiregularly on the set
of non-trivial conjugacy classes of~$A$. So any set $C\in\cla_G(A)$ other than $\{1\}$
is a disjoint union of exactly $k=|K|$ of these classes. This shows that
\qtnl{300509b}
|\cla_G(A)|=1+(|\cla(A)|-1)/k
\eqtn
where $\cla(A)$ is the set of conjugacy classes of~$A$. To complete the proof we note
that from (\ref{280509a}) it follows that the scheme $S$ is commutative if and only
if $\CC S=\pi(Z(\CC G))$, or equivalently if
$$
1+(|A|-1)/k=|S|=\dim(\CC S)=\dim(Z(\CC S))=\dim(\pi(Z(\CC G))).
$$
By (\ref{300509a}) and (\ref{300509b}) this is true if and only if $|\cla(A)|=|A|$,
i.e. the group $A$ is commutative.\bull

We would like to have a characterization of a pseudocyclic scheme in terms of its
intersection number array. For commutative case it was done in~\cite[Proposition~2.2.7]{BCN}.

\thrml{pscyclic}
The following two statements are equivalent:
\nmrt
\tm{1} $(\Omega,S)$ is a pseudocyclic scheme with $m_\SP/n_\SP=k$ for all $P\in\P^\#$,
\tm{2} $(\Omega,S)$ is an equivalenced scheme of valency $k$ with $c(r)=k-1$ for
all $r\in S^\#$.
\enmrt
Moreover, any pseudocyclic scheme with pairwise equal non-principal dimensions
of irreducible representations is commutative.
\ethrm
\proof Set $n=|\Omega|$ and $r=|S|$. We will use the following identity proved in
Proposition~3.4 and Lemma~3.8(i) of~\cite{AFM}:
\qtnl{eqafm1}
\sum_{s\in S}\frac{\reg(s^*)}{n_s}A(s)=n\sum_{P\in\P}\frac{n_P}{m_P}P.
\eqtn
where $\reg(s^*)=\sum_{t\in S}c_{s^*t}^t$. Suppose that $m_\SP/n_\SP=k$ for all
$P\in\P^\#$. Then taking into account that $I=J/n+\sum_{P\in\P^\#}P$ where $I$
is the identity matrix, from (\ref{eqafm1}) we obtain
\qtnl{eqafm}
\sum_{s\in S} \frac{\reg(s^\star)}{n_s} A(s) =
\frac{n}{k} I + \frac{k-1}{k}J.
\eqtn
By comparison of the diagonal and non-diagonal entries of the matrices in both sides of
this equality we conclude (after rearrangements) that
\qtnl{300509d}
k=\frac{n-1}{r-1},\qquad \reg(s)=\frac{n_s(k-1)}{k},\quad s\in S^\#.
\eqtn
Thus
$\reg(s)=n_s(n-r)/(n-1)$ for each $s\in S^\#$. Denote by $f$ be the great common divisor
of $n_s$, $s\in S^\#$. Since $f=\sum_{s\in S^\#}x_s n_s$
for some $x_s\in\ZZ$, we obtain $f(n-r)/(n-1)=\sum_{s\in S^\#} x_s\reg(s)\in\ZZ$.
Therefore $(r-1)f\ge n-1$, and hence
$$
n-1=\sum_{s\in S^\#}n_s\ge \sum_{s\in S^\#}f=(r-1)f\ge n-1.
$$
This implies that $n_s=f=(n-1)/(r-1)$ for all $s\in S^\#$. Due to (\ref{300509d}) this
means that $(\Omega,S)$ is an equivalenced scheme of valency $k$ with $c(r)=\reg(r)=k-1$ for
all $r\in S^\#$.\medskip

Assume now that $(\Omega,S)$ is an equivalenced scheme of valency $k$ with $c(r)=k-1$
for all $r\in S^\#$. Then $\reg(s)=c(r)$ for all~$r$. So the left side of~\eqref{eqafm1}
is a linear combination of the matrices $I$ and $J$. After multiplying both sides of that
equality by $P\in\P^\#$ we see that $m_P/n_P=k$. Thus the scheme $(\Omega,S)$ is
pseudocyclic with $m_P/n_P=k$ for all $P\in\P^\#$.\medskip

Let now $(\Omega,S)$ be a pseudocyclic scheme such that $n_P$ does not depend on $P\in\P^\#$.
Denote this number by~$a$. Then by the first part of the proof the scheme is
equivalenced of valency $k$ where $ka=m_P$ for each $P\in\P^\#$. So from
(\ref{190309a}) it follows that $|\P^\#|a^2=r-1=(n-1)/k$. Therefore $a$ is coprime to~$n$.
Taking into account that $n_s=k$ for all $s\in S^\#$, the Frame number of the scheme
$(\Omega,S)$ can be computed as follows
$$
n^r\frac{\prod_{s\in S}n_s}{\prod_{P\in\P}m_\SP^{n_\SP^2}}=
\frac{n^rk^{r-1}}{k^{r-1}a^{r-1}}=\frac{n^r}{a^{r-1}}.
$$
Since this number is an integer and $r>1$, we conclude that $a=1$. Thus
$n_P=1$ for all~$P$ and the scheme is commutative.\bull

There is a lot of equivalenced schemes $(\Omega,S)$ for which the group of algebraic
isomorphisms acts transitively on $S^\#$. These schemes were first studied by Ikuta, Ito
and Munemasa~\cite{IIM}, and include the cyclotomic schemes over finite
fields and the schemes of affine planes (see Section~\ref{040209a}). The following
statement shows that all of them are pseudocyclic.

\crllrl{040209b}
Let $(\Omega,S)$ be an equivalenced scheme. Suppose that a group of its algebraic
isomorphisms acts transitively on $S^\#$. Then $(\Omega,S)$ is a pseudocyclic scheme.
\ecrllr
\proof From the hypothesis it follows that the number $c(s)$ does not depend on $s\in S^\#$.
So by Lemma~\ref{180309u} this number equals $k-1$ and we are done by Theorem~\ref{pscyclic}.\bull

Sometimes one can construct a new pseudocyclic scheme by means of an appropriate
{\it algebraic fusion} defined as follows. Let $G$ be a group of algebraic isomorphisms
of a coherent configuration $(\Omega,S)$. Set
$$
S^G=\{s^G:\ s\in S\}
$$
where $s^G$ is the union of the relations $s^g$, $g\in G$. It is easily seen that
the pair $(\Omega,S^G)$ is a coherent configuration. Moreover, if the group $G$
is half-transitive on $S^\#$ and the coherent configuration $(\Omega,S)$ is
equivalenced, then $(\Omega,S^G)$ is an equivalenced scheme. An analog of this
statement holds for commutative pseudocyclic schemes.\footnote{Using Theorem~\ref{060209a}
enables us to reduce substantially the proofs in \cite[Section~3]{HX06}.}

\thrml{060209a}
Let $(\Omega,S)$ be a commutative pseudocyclic scheme of valency $k$ and $G$ a group of
algebraic isomorphisms of it. Suppose that $G$ acts semiregularly on $S^\#$. Then
$(\Omega,S^G)$ is a commutative pseudocyclic scheme of valency $km$ where $m=|G|$.
\ethrm
\proof We observe that $(\Omega,S^G)$ being a fusion of a commutative scheme is also
commutative. Since it is equivalenced of valency $km$ (see above), we have
\qtnl{010409z}
|\P^G|=|S^G|=1+(|S|-1)/m
\eqtn
where $\P^G$ is the set of central primitive idempotents of the algebra $\CC S^G$.
On the other hand, the group $G$ naturally acts as an automorphism group of
the algebra $\CC S$. Therefore the set $\P$ of its central primitive idempotents
is $G$-invariant. For $P\in\P$ denote by $P^G$ the sum of all $Q\in\P$ belonging
to the $G$-orbit containing~$P$. Then obviously the set $\P'=\{P^G:\ P\in\P\}$
consists of pairwise orthogonal central idempotents of the algebra $\CC S^G$.
This implies that
\qtnl{010409y}
|\P'|\le |\P^G|.
\eqtn
Finally, any $G$-orbit in $\P$ is of cardinality at most~$m$. Since  $P_0$
leaved fixed under $G$ and $|\P|=|S|$, this implies that
$$
|\P'|\ge 1+(|\P|-1)/m=1+(|S|-1)/m
$$
and the equality holds exactly when any $G$-orbit in $\P^\#$ is of size $m$. Together
with (\ref{010409z}) and (\ref{010409y}) this shows that $|\P^G|=|\P'|$. Therefore
$m_\SPG=mm_P=mk$ and $n_\SPG=1$ for all $P\in\P^\#$. Thus the scheme $(\Omega,S^G)$ is
pseudocyclic.\bull

{A schurian equivalenced non-regular scheme is nothing but the scheme of $3/2$-transitive group.}
In general, the latter is not a Frobenius
group. However, the following statement holds.

\thrml{020209d}
A schurian pseudocyclic scheme of valency $k>1$ and rank greater than $2(k-1)$
is a Frobenius scheme the automorphism group of which is a Frobenius group.
\ethrm
\proof Let $(\Omega,S)$ be a schurian pseudocyclic scheme of valency $k>1$ and
$G=\aut(\Omega,S)$. Then the set $\fix(g)$ of points left fixed by a nonidentity
permutation $g\in G$ does not coincide with~$\Omega$. So there exists $\alpha\in\Omega$
such that $\alpha^g\ne\alpha$. Since $r(\alpha,\beta)=r(\alpha^g,\beta)$
for all $\beta\in\fix(g)$, Theorem~\ref{pscyclic} implies that
$$
|\fix(g)|\le |\{\beta\in\Omega:\ r(\alpha,\beta)=r(\alpha^g,\beta)\}|=c(s)=k-1
$$
where $s=r(\alpha,\alpha^g)$. Thus the permutation character of $G$
takes the value in the set $\{0,\ldots,k-1\}$ on all nonidentity elements. Therefore by
\cite[Prop.1]{C01} a point stabilizer $G_\alpha$ has at most $2(k-1)-1$ non-regular
orbits. If $|S|>2(k-1)$, then at least one orbit of $G_\alpha$ is regular.
This implies that so are all non-trivial orbits. Thus $G$ is a Frobenius group.\bull

\section{Known examples of pseudocyclic schemes}\label{040209a}

\sbsnt{Schemes of rank~$3$.} Any pseudocyclic scheme of rank~$3$ arises from either
a conference matrix (symmetric case) or from a skew Hadamard matrix (antisymmetric case)
\cite{BCN,Pa92}. In any case the intersection number array is uniquely determined
by the degree of the scheme. In particular, two such schemes are algebraically isomorphic
if and only if they have the same degree. Since there is an infinite number of $n$
for which there are at least two non-equivalent conference matrices or non-equivalent
skew Hadamard matrices of order $n$, in general pseudocyclic schemes of rank~$3$ are not separable.
Similarly, one can see that most of them are non-schurian. For instance, it follows from \cite{B72}
that an antisymmetric pseudocyclic scheme of rank~$3$ is schurian if and only if
it is a scheme of a Paley tournament. Moreover, in \cite[p.75]{FKM} one can find
an infinite family of non-schurian pseudocyclic schemes of rank~$3$ satisfying the
$4$-condition.

\sbsnt{The Hollman schemes \cite[p.390]{BCN}.} Let $q>4$ be a power of $2$. Denote by
$\Omega$ the set of cyclic groups of order $q+1$ in the group $\PSL(2,q)$. The latter acts
transitively on $\Omega$ by conjugation and hence produces the scheme of degree
$(q^2-q)/2$. One can prove that this scheme is symmetric and pseudocyclic of valency
$q+1$. Some algebraic fusions of the Hollman scheme that are also pseudocyclic were
studied in~\cite{HX06}.

\sbsnt{The Passman schemes \cite{Pa67}.} Let $q$ be an
odd prime power and $G$ the group consisting the transformations
\qtnl{110209a}
\begin{pmatrix}
x \\
y \\
\end{pmatrix}
\rightarrow
\begin{pmatrix}
a & 0 \\
0 & \pm a^{-1}\\
\end{pmatrix}
\begin{pmatrix}
x \\
y \\
\end{pmatrix}
+
\begin{pmatrix}
b \\
c \\
\end{pmatrix},
\qquad
\begin{pmatrix}
x \\
y \\
\end{pmatrix}
\rightarrow
\begin{pmatrix}
0 & a \\
\pm a^{-1} & 0\\
\end{pmatrix}
\begin{pmatrix}
x \\
y \\
\end{pmatrix}
+
\begin{pmatrix}
b \\
c \\
\end{pmatrix}
\eqtn
where $a,b,c\in\GF(q)$, and $a\ne 0$. Then $G$ is a $3/2$-transitive group acting on a
$2$-dimensional space over $\GF(q)$. The scheme of this group is equivalenced
of degree $q^2$ and valency $2(q-1)$. In fact, the Passman scheme is the algebraic fusion
of the Frobenius scheme of valency $q-1$ corresponding to the subgroup of
$G$ of order $(q-1)q^2$ consisting of the first family of permutations from~(\ref{110209a}).
Thus by Theorems~\ref{190309g} and~\ref{060209a} the Passman scheme is pseudocyclic.

\sbsnt{Cyclotomic schemes.}
Let $R$ be a finite local commutative ring with identity. Then its multiplicative
group $R^\times$ is the direct product of the Teichm\"uller group $\T$ and the group of
principal units~\cite{MD74}. The Teichm\"uller group is isomorphic to the multiplicative
group of the residue field of~$R$, and acts as a fixed point free automorphism
group of the additive group $R^+$ of~$R$. {Therefore for a given $K\le\T$,
 the scheme
$\cyc(K,R)$ of the group\\ $(R^+\rtimes K,R^+)$}
is a Frobenius scheme and hence pseudocyclic by Theorem~\ref{190309g}.
This example is a special case of a cyclotomic scheme over a finite commutative
ring~\cite{EP09}. In almost the same way one can construct a class of pseudocyclic
schemes where the ground ring $R$ is replaced by a near-field~\cite{BPR}, or even a
near-ring.

\sbsnt{Affine schemes.}\label{100209d} Let $\Omega$ be a point set of a finite affine
space $\A$ (see~\cite{BC95}). Denote by $S$ the partition of $\Omega\times\Omega$
containing $1_\Omega$ and such that two pairs
$(\alpha,\beta),(\alpha',\beta')\in\Omega\times\Omega$, $\alpha\ne\beta$,
$\alpha'\ne\beta'$, belong to the same class if and only if
the lines $\alpha\beta$ and $\alpha'\beta'$ are equal or parallel.
Then the pair $(\Omega,S)$ is a symmetric scheme and nonzero intersection
numbers $c_{rs}^t$ with $1_\Omega\not\in\{r,s\}$ are as follows:
\qtnl{070209a}
c_{rs}^t=\css
q-1,&\text{$r=s,\ t=1_\Omega$},\\
q-2,&\text{$r=s=t$},\\
1, &\text{$r\ne s,\ t\in rs$},\\
\ecss
\eqtn
where $q$ is the size of a line in $\A$ (called the order of~$\A$). It follows that
it is a pseudocyclic scheme of valency $q-1$. Each relation of an affine scheme is
an involution in a sense of \cite{Zi1}. It was shown in \cite{D00} that a scheme whose relations are involutions is an affine scheme.
Thus there is one-to-one correspondence between affine schemes and affine spaces.
It is straightforward
to prove that the schemes of affine spaces are isomorphic (resp. algebraically isomorphic)
if and only if the affine spaces are isomorphic (resp. have the same order).

\thrml{090209z}
The scheme of a finite affine space $\A$ is schurian if and only if $\A$ is Desarguesian.
\ethrm
\proof By the Veblen-Young theorem (cf.~\cite{BC95}) a finite affine space $\A$ is either
a non-Desarguesian affine plane, or the $n$-dimensional affine geometry $\AG(n,q)$
over $\GF(q)$. In the latter case the scheme of $\A$ coincides with the scheme of the
group
$$
TC\le\AGL(n,q)
$$
where $T$ is the translation group and $C$ is the centre of $\GL(n,q)$. This
proves the sufficiency part of the theorem. To prove the necessity assume that the
scheme of $\A$ is schurian. Then it satisfies the $4$-condition and the required
statement immediately follows from the lemma below.

\lmml{100209c}
Suppose that the scheme of an affine space $\A$ satisfies the $4$-condition. Then $\A$
is Desarguesian.
\elmm
\proof Let $(\Omega,S)$ be the scheme of $\A$. It suffices to verify that given seven
distinct points $\alpha,\alpha'$, $\beta,\beta'$, $\gamma,\gamma'$ and $\delta$, such that
$\alpha\alpha'$, $\beta\beta'$, and $\gamma\gamma'$ are distinct lines through $\delta$
and $\alpha\gamma$ is parallel to $\alpha'\gamma'$ and $\beta\gamma$ is parallel to
$\beta'\gamma'$, then $\alpha\beta$ is parallel to $\alpha'\beta'$ (see Fig~\ref{f8}).
\begin{figure}[h]
$\xymatrix@R=10pt@C=10pt@M=0pt@L=5pt{
& & & & & & \VRT{\delta} \ar@{-}[ddll] \ar@{-}[ddd] \ar@{-}[ddrr] & & & & & &\\
& & & & & & & & & & & &\\
& & & & \VRT{\alpha}\ar@{-}[drr] \ar@{-}[ddll]\ar@{-}[rrrr] & & & & \VRT{\gamma} \ar@{-}[ddrr]\ar@{-}[dll] & & & &\\
& & & & & & \VRT{\beta} \ar@{-}[dddd] & & & & & &\\
& & \VRT{\alpha'}\ar@{-}[rrrrrrrr] \ar@{-}[dddrrrr]& & & & & & & & \VRT{\gamma'}\ar@{-}[dddllll] & & \\
& & & & & & & & & & & & & &\\
& & & & & & & & & & & & & &\\
& & & & & & \VRT{\beta'} & & & & & & \\
}$
\caption{}\label{f8}
\end{figure}
However, since $\delta\gamma=\delta\gamma'$, we have $r(\delta,\gamma)=r(\delta,\gamma')$.
Due to the $4$-condition there exist points $\alpha'',\beta''$ such that the $4$-sets
$\Delta=\{\delta,\gamma,\alpha,\beta\}$ and $\Delta'=\{\delta,\gamma',\alpha'',\beta''\}$
have the same type with respect to the pairs $(\delta,\gamma)$ and $(\delta,\gamma')$
respectively. So $r(\alpha,\delta)=r(\alpha'',\delta)$ and
$r(\beta,\delta)=r(\beta'',\delta)$. This implies that
\qtnl{020409a}
\alpha''\in\alpha\delta=\alpha'\delta,\quad
\beta''\in\beta\delta=\beta'\delta.
\eqtn
On the other hand, $r(\alpha,\gamma)=r(\alpha'',\gamma')$ and
$r(\beta,\gamma)=r(\beta'',\gamma')$. Therefore $\alpha\gamma$ is parallel to
$\alpha''\gamma'$, and $\beta\gamma$ is parallel to $\beta''\gamma'$. Thus from
(\ref{020409a}) we conclude that $\alpha''=\alpha'$ and $\beta''=\beta'$. Since
also $r(\alpha,\beta)=r(\alpha'',\beta'')$, the line $\alpha\beta$ is parallel
to $\alpha''\beta''=\alpha'\beta'$, and we are done.\bull

\sbsnt{Amorphic schemes.}
A scheme $(\Omega,S)$ is called amorphic \cite{GIK} if any its fusion is a
scheme. It was shown in \cite{GIK} that all basis graphs of an amorphic scheme of rank at least
four are strongly regular either of Latin square type or of negative Latin square type.
If an amorphic scheme is equivalenced, then its group of algebraic automorphisms
is $\sym(S^\#)$. This implies (Corollary~\ref{040209b}) that $(\Omega,S)$ is pseudocyclic.
A scheme of an affine plane of order $q$ is an amorphic $(q-1)$-valenced scheme of rank $q+2$.
This yields us
the following statement.

\thrm
Let $q$ be the order of an affine plane. Then given a divisor $m$ of $q+1$ and a partition
of $\{1,\ldots,q+1\}$ in $m$ classes of cardinality $(q+1)/m$, there exists an
amorphic pseudocyclic scheme of degree $q^2$, valency $(q^2-1)/m$ and rank $m+1$.\bull
\ethrm

\section{One point extension of an equivalenced scheme.}\label{090209g}

\sbsnt{Splitting sets.}
Let $(\Omega,S)$ be an equivalenced scheme of valency~$k$. For each
(possibly equal) basis relations $u,v\in S^\#$ we define the {\it splitting set} of
them as follows:
\qtnl{150309a}
D(u,v)=\{w\in S^\#:\ (uu^*\,vv^*)\cap ww^*=\{1_\Omega\}\}.
\eqtn
It is easily seen that $D(u,v)=D(v,u)$ and $uu^*\cap ww^*=vv^*\cap ww^*=\{1_\Omega\}$ for
all $w\in D(u,v)$. Therefore from Lemma~\ref{141208j} it follows that
\qtnl{150309d}
c_{u^*w}^s\le 1\quad\text{and}\quad c_{w^*v}^s\le 1
\eqtn
for all $s\in S$. In particular, $|u^*w|=|w^*v|=k$.

\thrml{150309e}
Given $w\in D(u,v)$ the following statements hold:
\nmrt
\tm{1} $|u^*v|=k\ \Leftrightarrow\ u\in D(v,w)\ \Leftrightarrow\ v\in D(w,u)$,
\tm{2} $|ab\cap u^*v|=1$ for all $a\in u^*w$ and $b\in w^*v$.
\enmrt
\ethrm
\proof To prove statement (1) suppose that $|u^*v|=k$. Then $c_{u^*v}^t\le 1$ for
all~$t\in S$. So from Lemma~\ref{141208j} it follows that $uu^*\cap vv^*=\{1_\Omega\}$.
On the other hand, given a relation $s\in (vv^*\,ww^*)\cap uu^*$ one can find points
$\alpha,\beta\in\Omega$ to have the configuration at Fig.~\ref{f1}.
\begin{figure}[h]
\grphp{
& & \EVR \ar[ddll]_*{u} \ar[dddd]^*{s} \ar[rrr]^*{v} & & & \EVR{} & &\\
& & & & & & & \\
\EVR{} & & & & & & & \VRT{\beta} \ar[lldd]^*{w} \ar[lluu]^*{v}\\
& & & & & & & \\
& & \VRT{\alpha} \ar[lluu]^*{u} \ar[rrr]^*{w} & & & \EVR{} & & \\
}
\caption{}\label{f1}
\end{figure}
So $r(\alpha,\beta)\in (uu^*\,vv^*)\cap ww^*=\{1_\Omega\}$ whence it follows that
$\alpha=\beta$. Therefore $s\in vv^*\cap uu^*=\{1_\Omega\}$. Thus $s=1_\Omega$, and
hence $u\in D(v,w)$. Conversely, if $u\in D(v,w)$, then $(vv^*\,ww^*)\cap uu^*=\{1_\Omega\}$.
So $vv^*\cap uu^*=\{1_\Omega\}$, and hence $|u^*v|=k$ by Lemma~\ref{141208j}.
This completes the proof of the first equivalence in statement~(1). The second
equivalence immediately follows from the first one by interchanging $u$ and $v$
because $D(u,v)=D(v,u)$.\medskip

To prove statement (2) let $a\in u^*w$ and $b\in w^*v$. Then $w\in ua\cap vb^*$. This
implies that $|ua\cap vb^*|\ge 1$, and hence $|ab\cap u^*v|\ge 1$. Thus it suffices to
verify that $|ab\cap u^*v|\le 1$. To do this we need the following auxiliary statement.

\lmml{150309f}
Given a relation $s\in ab\cap u^*v$ and points $\alpha,\beta,\gamma\in\Omega$ such that
$r(\alpha,\beta)=s$, $r(\alpha,\gamma)=a$ and $r(\gamma,\beta)=b$ there exists
a unique point $\delta\in\Omega$ for which $r(\delta,\alpha)=u$, $r(\delta,\beta)=v$
and $r(\delta,\gamma)=w$ (see Fig.~\ref{f2})
\elmm
\begin{figure}[h]
\grphp{
& & & \VRT{\alpha} \ar[ddddlll]_*{s} \ar[ddddrrr]^*{a} & & & \\
& & & & & & \\
& & & & & & \\
& & & \VRT{\delta} \ar[uuu]_*{u} \ar[llld]_*{v} \ar[rrrd]^*{w} & & & \\
\VRT{\beta} & & & & & & \VRT{\gamma} \ar[llllll]^*{b} \\
}
\caption{}\label{f2}
\end{figure}

\proof Since $a\in u^*w$ and $b\in w^*v$, there exist points $\lambda$, $\mu$, $\nu$
such that $r(\lambda,\alpha)=u$, $r(\lambda,\beta)=v$, $r(\mu,\alpha)=u$, $r(\mu,\gamma)=w$
and $r(\nu,\gamma)=w$, $r(\nu,\beta)=v$ (see Fig.~\ref{f3}).
\begin{figure}[h]
\grphp{
\VRT{\lambda} \ar[rrrr]^*{u} \ar[ddrr]_*{v} & & & & \VRT{\alpha} \ar[ddll]_*{s} \ar[ddrr]^*{a} & & & & \VRT{\mu} \ar[llll]_*{u} \ar[ddll]^*{w}\\
& & & & & & & &\\
& & \VRT{\beta} & & & & \VRT{\gamma} \ar[llll]^*{b}\\
& & & & & & & &\\
& & & & \VRT{\nu} \ar[lluu]^*{v} \ar[rruu]_*{w} & & & & \\
}
\caption{}\label{f3}
\end{figure}
Now $r(\mu,\nu)\in (uu^*\,vv^*)\cap ww^*=\{1_\Omega\}$. Thus $\mu=\nu$. Denote this
point by $\delta$. Then for $\delta$ the statement of the lemma holds. To prove the
uniqueness we note if $\delta_1$ and $\delta_2$ are two points forming Fig.~\ref{f2},
then the relation $r(\delta_1,\delta_2)$ belongs to the set $ww^*\cap uu^*\cap vv^*=\{1_\Omega\}$, and hence
$\delta_1=\delta_2$.\bull

To complete the proof of Theorem~\ref{150309e} suppose that $ab\cap u^*v\supset\{s_1,s_2\}$
with $s_1\ne s_2$. Then there exist points $\alpha,\gamma,\beta_1,\beta_2\in\Omega$,
such that $\beta_1\ne\beta_2$, $r(\alpha,\gamma)=a$ and $r(\gamma,\beta_i)=b$,
$r(\alpha,\beta_i)=s_i$ for $i=1,2$. By Lemma~\ref{150309f} with $s=s_i$ one can
find a point $\delta_i$ for which $r(\delta_i,\alpha)=u$, $r(\delta_i,\beta)=v$
and $r(\delta_i,\gamma)=w$ (see Fig.~\ref{f4}).
\begin{figure}[h]
\grphp{
& & & & & & \VRT{\gamma} \ar[dddllllll]_*{b} \ar[dddrrrrrr]^*{b} & & & & & & \\
& & & & & & & & & & & & \\
& & & & & & & & & & & &\\
\VRT{\beta_1} & & & & \VRT{\delta_1} \ar[llll]^*{v} \ar[uuurr]_*{w} \ar[dddrr]^*{u} & & & & \VRT{\delta_2} \ar[rrrr]_*{v} \ar[uuull]^*{w} \ar[dddll]_*{u} & & & &\VRT{\beta_2}\\
& & & & & & & & & & & & \\
& & & & & & & & & & & & \\
& & & & & & \VRT{\alpha} \ar[uuuuuu]_*{a} \ar[lllllluuu]^*{s_1} \ar[rrrrrruuu]_*{s_2}\\
}
\caption{}\label{f4}
\end{figure}
Since the relation $r(\delta_1,\delta_2)$ belongs to the set $ww^*\cap uu^*=\{1_\Omega\}$,
we have $\delta_1=\delta_2$. Denote this point by $\delta$. Then $r(\delta,\beta_i)=v$.
Since $r(\beta_i,\gamma)=b^*$, $r(\delta,\gamma)=w$ and $\beta_1\ne\beta_2$, this
implies that $c_{vb^*}^w\ge 2$. So by equivalencity $(X,\Omega)$ and (\ref{150309c})
we conclude that $c_{w^*v}^b\ge 2$ which contradicts to~(\ref{150309d}).\bull

\sbsnt{A one point extension.} By means of splitting sets we are going to find
explicitly the $\alpha_0$-extension of the scheme $(\Omega,S)$ for any point $\alpha_0\in\Omega$.
To do this for any relations $u,v\in S$ set
\qtnl{210109f}
S(u,v)=S_{\alpha_0}(u,v)=\{s\cap(\alpha_0u\times\alpha_0v):\ s\in u^*v\}.
\eqtn
It is easily seen that the union of relations from $S(u,v)$ coincides with the set
$\alpha_0u\times\alpha_0v$. Suppose that $w\in D(u,v)$. Then from statement~(2) of
Theorem~\ref{150309e} it follows that
\qtnl{180309c}
S(u,v)\subset S(u,v;w)^\cup
\eqtn
where $S(u,v;w)=S_{\alpha_0}(u,v;w)=S(u,w)\cdot S(w,v)$. Suppose, in addition, that $|u^*v|=k$.
Then any relation in $S(u,v)$ is of cardinality $k$. Since the same is true for the
relations in $S(u,v;w)$, we conclude that
\qtnl{040409a}
(w\in D(u,v)\quad \&\quad |u^*v|=k)\ \Rightarrow\ S(u,v)=S(u,v;w).
\eqtn
In particular, in this case the relations from $S(u,v;w)$ form a
partition of the set $\alpha_0u\times\alpha_0v$. For arbitrary $u$ and $v$ we will prove
the latter only under the following additional assumption:
\qtnl{180309a}
\bigcap_{\substack{a\in\{u,v\}\\b\in\{w,w'\}}}D(a,b)\ne\emptyset,\qquad w,w'\in D(u,v).
\eqtn
In this case obviously $D(u,v)\ne\emptyset$.

\lmml{260109a}
Let $u,v\in S^\#$ be such that (\ref{180309a}) holds. Then the set $S(u,v;w)$ forms a
partition of the set $\alpha_0u\times\alpha_0v$, and this partition does not depend on the choice of
the relation $w\in D(u,v)$.
\elmm
\proof Let $w,w'\in D(u,v)$. Then $|u^*w|=|u^*w'|=|v^*w|=|v^*w'|=k$. On the other
hand, due to (\ref{180309a}) one can find a relation
$$
t\in D(u,w)\cap D(u,w')\cap D(v,w)\cap D(v,w').
$$
So from (\ref{040409a}) it follows that
\qtnl{170309c}
S(x,y)=S(x,y;t)=S(x,t)\cdot S(t,y),\qquad x\in\{u,v\},\ y\in\{w,w'\}.
\eqtn
Let $a\in S(u,t)$ and $b\in S(t,v)$. Since obviously $b^*\in S(v,t)$ by~(\ref{170309c})
with $(x,y)=(u,w)$ and $(x,y)=(v,w)$ we obtain
$$
a\cdot c\in S(u,w),\quad b^*\cdot c\in S(v,w),\qquad c\in S(t,w).
$$
So the element $a\cdot b\in S(u,v;w)$ has at least $k$ different representations (one for each
choice of $c$) of the form $a\cdot b=(a\cdot c)(c^*\cdot b)$ with $c^*\cdot b\in S(w,v)$. Since
$|S(u,w)|=|S(w,v)|=k$ this implies that $|S(u,v;w)|=k$. Thus $S(u,v;w)$ is a partition
of $\alpha_0u\times \alpha_0v$ and $S(u,v;w)=S(u,v;t)$. Similarly, using
equalities~(\ref{170309c}) with $(x,y)=(u,w')$ and $(x,y)=(v,w')$ one can prove that
$S(u,v;w')$ is a partition of $\alpha_0u\times\alpha_0v$ and that $S(u,v;w')=S(u,v;t)$.
Thus
$$
S(u,v;w)=S(u,v;t)=S(u,v;w')
$$
and we are done.\bull

By Lemma~\ref{260109a} we can define a uniquely determined partition of the set
$\Omega_0\times\Omega_0$ where $\Omega_0=\Omega\setminus\{\alpha_0\}$, as follows
\qtnl{270109a}
S_0=S_0(\alpha_0)=\bigcup_{w\in D(u,v)}S_{\alpha_0}(u,v;w).
\eqtn
It is easily seen that $1_{\alpha_0u}\in S_0$ for all $u\in S^\#$. Therefore
$1_{\Omega_0}\in S_0^\cup$. Besides, since obviously $S(u,v;w)^*=S(v,u;w)$
for all $u,v,w$, the partition $S_0$ is closed with respect to $*$. Finally,
the condition (\ref{040409b}) is satisfied for $S=S_0$. Therefore if
$(\Omega_0,S_0)$ is a coherent configuration, then it is semiregular.
We will prove the former under the following assumption:
\qtnl{210309a}
D(u,v)\cap D(v,w)\cap D(w,u)\ne\emptyset
\eqtn
for all $u,v,w\in S^\#$. Below we set
$$
S_1(\alpha_0)=\{\{\alpha_0\}\times\alpha_0v: v\in S\}\cup \{\alpha_0v \times\{\alpha_0\}: v\in S\}.
$$

\thrml{280109c}
Let $(\Omega,S)$ be an equivalenced scheme satisfying conditions~(\ref{180309a})
and~(\ref{210309a}) for all $u,v\in S^\#$. Then
$$
S_{\alpha_0}=S_0(\alpha_0)\cup S_1(\alpha_0),\qquad \alpha_0\in\Omega.
$$
In particular, the $\alpha_0$-extension of $(\Omega,S)$ is a semiregular coherent
configuration on $\Omega_0$ the fibers of which are $\alpha_0u$, $u\in S$.
\ethrm
\proof It suffices to verify that $(\Omega_0,S_0)$ is a coherent configuration. Indeed,
if it is so and $S'=S_0(\alpha_0)\cup S_1(\alpha_0)$, then obviously $(\Omega,S')$ is a
semiregular coherent configuration the fibers of which are $\alpha_0u$, $u\in S$.
Therefore due to~(\ref{180309c}) we have $1_{\alpha_0}\in S'$ and $S\subset (S')^\cup$.
By the minimality of the $\alpha_0$-extension this implies that
\qtnl{210309c}
S_{\alpha_0}\subset (S')^\cup.
\eqtn
On the other hand, it is easily seen that the set $S(u,v)$, and hence the set
$S(u,v;w)=S(u,w)\cdot S(w,v)$ is contained in $(S_{\alpha_0})^\cup$ for all $u,v,w\in S$. Therefore
$S'\subset (S_{\alpha_0})^\cup$. Together with (\ref{210309c}) this shows that
$(S_{\alpha_0})^\cup=(S')^\cup$. Thus $S_{\alpha_0}=S'$ and we are done.\medskip

Let us prove that $(\Omega_0,S_0)$ is a coherent configuration. We observe that due
to~(\ref{180309a}) Lemma~\ref{260109a} implies that $S_0$ is a partition of
$\Omega_0\times\Omega_0$. Thus it suffices to verify that if $b,c\in S_0$ with
$b\cdot c\ne\emptyset$, then $b\cdot c\in S_0$ (see the remarks after~(\ref{270109a})). However,
for such $b$ and $c$ we have
$$
b\subset \alpha_0u\times\alpha_0v,\quad
c\subset \alpha_0v\times\alpha_0w
$$
for appropriate $u,v,w\in S^\#$. By condition (\ref{210309a}) one can find
$$
t\in D(u,v)\cap D(v,w)\cap D(w,u).
$$
Since $b=a_1\cdot a_2$ for some $a_1\in S(u,v;t)$ and $a_2\in S(v,w;t)$, this implies that
$|S(v,t)|=|S(w,t)|=k$ where $k$ is the valency of $(\Omega,S)$. So by statement~(1) of
Theorem~\ref{150309e} we have $w\in D(v,t)$. Therefore $a_2^*\in S(v,t;w)$ and hence
there exists $a_3\in S(t,w)$ such that $c=a_2^*\cdot a_3$. Thus
$$
b\cdot c=(a^{}_1\cdot a^{}_2)\cdot(a_2^*\cdot a^{}_3)=
a^{}_1\cdot(a^{}_2\cdot a_2^*)\cdot a^{}_3\in S(u,w;t).
$$
This means that $b\cdot c\in S_0$ which completes the proof.\bul

\section{One point extension of an algebraic isomorphism}\label{310309b}

We keep the notation of Section~\ref{090209g}. Let $(\Omega',S')$
be a scheme and let
$$
\varphi:(\Omega,S)\to(\Omega',S'),\ s\mapsto s'
$$
be an algebraic isomorphism. Then obviously $(\Omega',S')$ is an equivalenced scheme
of valency $k$, and
$$
w\in D(u,v)\ \Leftrightarrow\ w'\in D(u',v'),\qquad u,v\in S^\#.
$$
Let us fix a point $\alpha'_0\in\Omega'$.
Take $u,v\in S$. Since $c_{u^{}v^{}}^{w^{}}=c_{u'v'}^{w'}$ for all $w\in S$, the
mapping $\varphi$ induces a bijection
$$
\varphi_{u,v}:S_{\alpha^{}_0}(u,v)\to S_{\alpha'_0}(u',v'),\quad a\mapsto a'
$$
where $a=s\cap(\alpha_0u\times\alpha_0v)$ and $a'=s'\cap (\alpha'_0u'\times\alpha'_0v')$.
Below we set $S'(u',v')=S_{\alpha'_0}(u',v')$.

\lmml{260109o}
Let $u,v\in S^\#$ be such that (\ref{180309a}) holds. Then for any relations
$w_1,w_2\in D(u,v)$ we have
$$
b_1\cdot c_1=b_2\cdot c_2\ \Rightarrow\ b'_1\cdot c'_1=b'_2\cdot c'_2
$$
for all $b_i\in S(u,w_i)$ and $c_i\in S(w_i,v)$, $i=1,2$.
\elmm
\proof Suppose that $b_1\cdot c_1=b_2\cdot c_2$ where $b_i\in S(u,w_i)$ and $c_i\in S(w_i,v)$, $i=1,2$.
Then $|u^*w_1|=|u^*w_2|=|v^*w_1|=|v^*w_2|=k$. Besides, by~(\ref{180309a}) one can find a
relation
$$
t\in D(u,w_1)\cap D(u,w_2)\cap D(v,w_1)\cap D(v,w_2).
$$
So from (\ref{040409a}) it follows that
\qtnl{290309g}
S(x,y;t)=S(x,t)\cdot S(t,y),\qquad x\in\{u,v\},\ y\in\{w_1,w_2\}.
\eqtn
Take any $a_1\in S(u,t)$. By (\ref{290309g}) we have $d_1:=a_1^*b_1\in S(t,w_1)$ and
$d_2:=a_1^*\cdot b_2\in S(t,w_2)$.  Since $w_1\in D(t,v)$, we also have
$a_2:=d_1\cdot c_1\in S(t,v)$. Thus
$$
a_1\cdot a_2=(a_1\cdot d_1)\cdot (d^*_1\cdot a_2)=b_1\cdot c_1=b_2\cdot c_2=a_1\cdot d_2\cdot c_2
$$
whence it follows that $c_2=d_2^*\cdot a_2$ (Figure~\ref{f5}).
\begin{figure}[h]
\grphp{
& & & \VRT{w_1} \ar[dddrrr]^*{c_1} & & & \\
& & & & & & \\
& & & & & & \\
\VRT{u}\ar[rrruuu]^*{b_1}\ar[rrr]_*{a_1}\ar[rrrddd]_*{b_2} & & &
\VRT{t}\ar[uuu]_*{d_1}\ar[rrr]_*{a_2}\ar[ddd]^*{d_2} & & &
\VRT{v} \\
& & & & & & \\
& & & & & &\\
& & & \VRT{w_2} \ar[uuurrr]_*{c_2} & & & \\
}
\caption{}\label{f5}
\end{figure}
By Theorem~\ref{150309e} this implies that
$$
b'_1\cdot c'_1=(a_1\cdot d_1)'\cdot (d^*_1\cdot a_2)'=
a'_1\cdot (d'_1\cdot (d'_1)^*)\cdot a'_2=a'_1\cdot a'_2=
$$
$$
(b_2\cdot d_2^*)'\cdot (d_2\cdot c_2)'=
b'_2\cdot (d'_2\cdot (d^*_2)')\cdot c'_2=b'_2\cdot c'_2
$$
and we are done.\bull

Let us define a mapping $\varphi_0:S^{}_0\to S'_0$ where $S_0=S^{}_{\alpha^{}_0}$ and
$S'_0=S'_{\alpha'_0}$, as follows. Take $s\in S_0$. Then $s\subset\alpha_0u\times\alpha_0v$
for some $u,v\in S$. Set
$$
s^{\varphi_0}=\css
s^{\varphi_{u,v}}, &\text{if $|u^*v|=k$ or $1_\Omega\in\{u,v\}$},\\
s_1^{\varphi_{u,w}}\cdot s_2^{\varphi_{w,v}}, &\text{otherwise},
\ecss
$$
where $w\in D(u,v)$ and $s_1\in S(u,w)$, $s_2\in S(w,v)$ are such that $s_1\cdot s_2=s$.
By Theorem~\ref{280109c} and Lemma~\ref{260109o} the mapping $\varphi_0$ is a correctly
defined bijection.

\thrml{280109e}
Let $(\Omega,S)$ be an equivalenced scheme satisfying conditions~(\ref{180309a})
and~(\ref{210309a}) for all $u,v\in S^\#$, and $\varphi:(\Omega,S)\to(\Omega',S')$ an
algebraic isomorphism. Then $\varphi_0$ is the $(\alpha^{}_0,\alpha'_0)$-extension
of $\varphi$.
\ethrm
\proof We observe that relations (\ref{290309z}) hold by statement~(2) of
Theorem~\ref{150309e}. Therefore it suffices to verify that
\qtnl{300309a}
(b\cdot c)^{\varphi_0}=b^{\varphi_0}\cdot c^{\varphi_0},\qquad b,c\in S_0,\ b\cdot c\ne\emptyset.
\eqtn
To do this take such relations $b$ and $c$. Then there exist $u,v,w\in S^\#$
such that
$$
b\subset \alpha_0u\times\alpha_0v,\quad
c\subset \alpha_0v\times\alpha_0w.
$$
By condition (\ref{210309a}) one can find $t\in D(u,v)\cap D(v,w)\cap D(w,u)$.
Since $b\in S_0$, there exist $a_1\in S(u,t)$ and $a_2\in S(t,v)$ such that
$b=a_1\cdot a_2$. However, $|S(u,t)|=|S(t,v)|=k$. Therefore
$b^{\varphi_0}=a_1^{\varphi_0}\cdot a_2^{\varphi_0}$ (if $|u^*v|=k$, then this follows
from statement~(2) of Theorem~\ref{150309e}; otherwise this immediately follows
from the definition of~$\varphi_0$). On the other hand, we have
$|S(v,t)|=|S(w,t)|=k$. Therefore $a_2^*\in S(v,t)$ and there exists $a_3\in S(t,w)$ such
that $c=a_2^*\cdot a_3$ (see Figure~\ref{f6}).
\begin{figure}[h]
\grphp{
& & & \VRT{s} \ar[ddddrrr]^*{a_3}\ar[dddd]^*{a_2} & & & \\
& & & & & & \\
& & & & & & \\
& & & & & & \\
\VRT{u}\ar[rrruuuu]^*{a_1}\ar[rrr]_*{b} & & &
\VRT{v}\ar[rrr]_*{c} & & &
\VRT{w} \\
}
\caption{}\label{f6}
\end{figure}
As above one can see that $c^{\varphi_0}=(a_2^*)^{\varphi_0}\cdot a_3^{\varphi_0}$
and that $(a_1\cdot a_3)^{\varphi_0}=(a_1)^{\varphi_0}\cdot a_3^{\varphi_0}$.
Thus
$$
(b\cdot c)^{\varphi_0}=(a_1\cdot a_2\cdot a_2^*\cdot a_3)^{\varphi_0}=(a_1\cdot a_3)^{\varphi_0}=
(a_1)^{\varphi_0}\cdot a_3^{\varphi_0}=
$$
$$
(a_1)^{\varphi_0}\cdot a_2^{\varphi_0}\cdot (a_2^*)^{\varphi_0}\cdot a_3^{\varphi_0}=
(a_1\cdot a_2)^{\varphi_0}\cdot (a_2^*\cdot a_3)^{\varphi_0}=b^{\varphi_0}\cdot c^{\varphi_0}
$$
which completes the proof.\bul

\section{Schurity and separability of equivalenced schemes}\label{310309c}

\sbsnt{Schurity.}\label{300309c} The conclusion of Theorem~\ref{280109c} gives a
sufficient condition for a scheme to be schurian.

\thrml{250309a}
Let $(\Omega,S)$ be a scheme such that for each $\alpha\in\Omega$ the coherent
configuration $(\Omega,S_\alpha)$ is semiregular on $\Omega\setminus\{\alpha\}$
and its fibers are $\alpha s$, $s\in S$. Then $(\Omega,S)$ is a regular or
Frobenius scheme.
\ethrm
\proof Let $\alpha\in\Omega$. Then by Theorem~\ref{210109e} the coherent configuration
$(\Omega,S_\alpha)$ is schurian, and hence due to (\ref{210109c}) any set $\alpha s$,
$s\in S$, is the orbit of the group $G_\alpha$ where $G=\aut(\Omega,S)$. Since obviously
$G_{\alpha,\beta}=\id_{\Omega}$ for all $\alpha\ne\beta$, it suffices to verify that
the group~$G$ is transitive. To do this we note that semiregularity of $G_\alpha$
on $\Omega\setminus\{\alpha\}$ implies that
\qtnl{250309h}
k:=|G_\alpha|=|\alpha s|,\qquad \alpha\in\Omega,\ s\in S^\#.
\eqtn
Let $\Delta\in\orb(G,\Omega)$. Since any orbit of $G$ is a disjoint union of some orbits
of the group $G_\alpha$, the number $|\Delta|$ is divided by $k$ if and only if
$\alpha\not\in\Delta$. However, this is impossible if $\Delta\ne\Omega$. Thus $G$
is transitive.\bull

Let $(\Omega,S)$ be the scheme of an affine space $\A$. Then from (\ref{070209a}) it
follows that $D(u,v)=S^\#\setminus uv$ for all $u,v\in S^\#$. Therefore
conditions~(\ref{180309a}) and~(\ref{210309a}) are satisfied whenever the dimension
of $\A$ is at least~$3$. Thus in this case by Theorems~\ref{280109c} and~\ref{250309a}
the scheme $(\Omega,S)$ is schurian. Since the scheme is imprimitive, its automorphism
group is a Frobenius group in its natural permutation representation. Using properties
of Frobenius groups one can prove, without using Veblen-Young Theorem, that the scheme
$(\Omega,S)$ is an affine scheme of a Desarguesian affine space.

\sbsnt{Separability.} Similarly to the previous subsection the conclusion of
Theorem~\ref{280109e} gives a sufficient condition for a scheme to be separable.

\thrml{290309q}
In the condition of Theorem~\ref{250309a} the scheme $(\Omega,S)$ is separable.
\ethrm
\proof From Theorem~\ref{280109e} it follows that any algebraic isomorphism
$\varphi:(\Omega,S)\to(\Omega,S')$ has the $(\alpha,\alpha')$-extension
$$
\varphi_{\alpha,\alpha'}:(\Omega,S^{}_{\alpha^{}})\to(\Omega',S'_{\alpha'})
$$
for all $(\alpha,\alpha')\in\Omega\times\Omega'$. By Theorem~\ref{280109c} the
coherent configuration $(\Omega,S_\alpha)$ is semiregular on $\Omega\setminus\{\alpha\}$,
and hence is separable (Theorem~\ref{210109e}). This implies that the algebraic
isomorphism $\varphi_{\alpha,\alpha'}$ is induced by an isomorphism. Due to
the right-hand side of~(\ref{290309z}) this shows that the same isomorphism
induces the algebraic isomorphism~$\varphi$. Thus any algebraic isomorphism of
the scheme $(\Omega,S)$ is induced by an isomorphism and hence this scheme is
separable.\bull

\sbsnt{Pseudocyclic schemes} Let $(\Omega,S)$ be an equivalenced scheme
of valency $k$ and indistinguishing number $c$.

\lmml{210309d}
$|\ov D(u,v)|< ck^3$, $u,v\in S^\#$ where $\ov D(u,v)=S^\#\setminus D(u,v)$.
\elmm
\proof Let $u,v\in S^\#$.  Then it is easily seen that $|uu^*\,vv^*|\le k^3$. Besides,
given $t\in uu^*\,vv^*$, $t\in S^\#$, there exist at most $c(t)-1\le c$ relations
$w\in S^\#$ such that $t\in ww^*$. Thus for at most $ck^3$ relations $w$ the set
$(uu^*\,vv^*)\cap ww^*$ contains an element $t\in S^\#$. This means that $|\ov D(u,v)|<ck^3$.\bull

Suppose that the scheme $(\Omega,S)$ is pseudocyclic. Then by Theorem~\ref{pscyclic}
we have $c=k-1$. If, in addition, $|S|\ge 4ck^3$, then by Lemma~\ref{210309d} conditions~(\ref{180309a})
and~(\ref{210309a}) are satisfied for all $u,v\in S^\#$. By Theorem~\ref{280109c}
this implies that given $\alpha\in\Omega$ the $\alpha$-extension of the scheme
$(\Omega,S)$ is a semiregular coherent configuration on $\Omega\setminus\{\alpha\}$ the
fibers of which are $\alpha u$, $u\in S$. So by Theorems~\ref{250309a} and~\ref{290309q}
we obtain the following statement.

\thrml{220309h}
Any pseudocyclic scheme of valency $k>1$ and rank at least $4(k-1)k^3$ is a separable
Frobenius scheme.\bull
\ethrm

\section{Miscellaneous}\label{310309d}

\subsection{$2$-designs}
Any commutative pseudocyclic scheme of of valency $k$ on $n$ points produces a
$2-(n,k,k-1)$-design \cite[Corollary~2.2.8]{BCN}. The same is also true in the
non-commutative case (Theorem~\ref{180309e}). It would be interesting to study
these designs in detail.

\thrml{180309e}
A scheme $(\Omega,S)$ on $n$ points is pseudocyclic of valency $k$ if and only if the
pair $(\Omega,B)$ with $B=\{\alpha s: \alpha\in\Omega,\ s\in S^\#\}$ is a
$2-(n,k,k-1)$-design.
\ethrm
\proof The pair $(\Omega,B)$ is an $2-(n,k,k-1)$-design if and only if $|\alpha s|=k$ for all
$\alpha\in\Omega$ and $s\in S^\#$, and the number of blocks $\alpha s\in B$ containing two
distinct points $\beta,\gamma\in\Omega$ coincides with $k-1$. However, the number of these
blocks is obviously equals $c(s)$ where $s=r(\beta,\gamma)$. Thus the required statement
follows from Theorem~\ref{pscyclic}.\bull

\sbsn
It was proved in \cite[Lemma 5.13]{EP99ce} that a one point extension of an imprimitive
equivalenced scheme is ``almost semiregular''. The following statement
shows that the imprimitivity condition can be removed for pseudocyclic schemes
the rank of which is much more than the valency.

\thrml{141208a}
Let $(\Omega,S)$ be a pseudocyclic scheme of valency $k$. Suppose that $|S|>2k(k-1)+2$.
Then given $\alpha\in\Omega$ the coherent configuration $(\Omega,S_\alpha)$ is semiregular
on $\Omega\setminus\{\alpha\}$.
\ethrm
\proof For a relation $u\in S^\#$ set $E(u)=\{v\in S^\#:\ |u^*v|=k\}$.
Since by Theorem~\ref{pscyclic} the indistinguishing number $c$ of the scheme $(\Omega,S)$ is
$k-1$, we have
$$
|S\setminus E(u)|=\sum_{b\in uu^*\setminus\{1_\Omega\}}|\{v\in S:\ b\in vv^*\}|\le
|uu^*|c\le k(k-1).
$$
By the theorem hypothesis this implies that
\qtnl{100209a}
|E(u)\cap E(v)|\ge |S|-2(k-1)k-2>0
\eqtn
for all non-equal $u,v\in S^\#$. Set
$S_0=\{s\in S_\alpha:\ s\subset \Omega_0\times\Omega_0\}$ where
$\Omega_0=\Omega\setminus\{\alpha\}$. To complete the
proof it suffices to verify that each $s_0\in S_0$ is contained in some $s\in S_0^\cup$
for which $|\beta s|\le 1$ for all $\beta\in\Omega$. However, it is easy to see that
given $s_0\in S_0$ one can find $u,v\in S^\#$ such that
$$
s_0\subset \alpha u\times\alpha v.
$$
Due to (\ref{100209a}) there exists $w\in E(u)\cap E(v)$. Since
$S_\alpha(u,w)$, $S_\alpha(w,v)\subset S_0^\cup$, we have
$S_\alpha(u,v;w)\subset S_0^\cup$. Thus $s_0\subset s$ for some $s\in S_\alpha(u,v;w)$.
To complete the proof it suffices to note that $|\beta s|\le 1$ for $\beta\in\Omega$.\bull

In the condition of Theorem~\ref{141208a} for any point
$\beta\in\Omega\setminus\{\alpha\}$ we have $|\beta s|\le 1$ for all $s\in S_\alpha$.
Therefore the scheme $(\Omega,S_\alpha)$ is $1$-regular
in the sense of~\cite{EP03}. Thus by Theorem~9.3 of this paper we obtain the following
statement.

\crllrl{100209b}
In the condition of Theorem~\ref{141208a} any one point extension of the scheme
$(\Omega,S)$ is schurian and separable.\bull
\ecrllr

We complete the subsection by making a remark that Corollary~\ref{100209b} together with
\cite[Theorem~4.6]{EP00} shows that any pseudocyclic scheme of valency $k$ and rank
$O(k^2)$ is $2$-schurian and $2$-separable in the sense of~\cite{EP09}.

\subsection{Affine schemes}
The following characterization of the affine schemes was known in commutative
case (see~\cite{D00}).

\thrm
Let $(\Omega,S)$ be a scheme with $n_s\geq 3$ for each $s\in S^\#$. Then it is
the scheme of an affine space if and only if $c_{rs}^t\le 1$ for all $r,s,t\in S$
such that $r\ne s^*$.
\ethrm
\proof The necessity follows from~(\ref{070209a}). To prove the sufficiency let
$s\in S^\#$. Then $|ss^*|>1$ because $n_s>1$, and $ss^*\cap rr^*=\{1_\Omega\}$ for all
$r\ne s$ (Lemma~\ref{141208j}). Thus
$$
|S^\#|\ge |\bigcup_{s\in S^\#}ss^*\setminus\{1_\Omega\}|=
\sum_{s\in S^\#}|ss^*\setminus\{1_\Omega\}|\ge |S^\#|.
$$
This implies that there exists a bijection $s\mapsto s'$ from $S^\#$ to itself such that
$ss^*=\{1_\Omega,s'\}$. In particular, the scheme $(\Omega,S)$ is symmetric and
$c_{s's}^s=n_s-1\geq 2$. Therefore $s'=s$ for each $s\in S^\#$. Thus each relation
from $S$ is an equivalence relation minus a diagonal. All such schemes were classified
in~\cite{D00} where it was proved that each of them is the scheme of an affine space.\bull

\subsection{QI-groups}
A permutation group $G\leq\sym(\Omega)$ is called {\it QI} \cite{cameron} if its
permutation module $\Q \Omega$ is a direct sum of a one dimensional module and an
irreducible one, say ${\Q \Omega}^0$. Notice that $(\Q \Omega)^0$ consists of all
vectors with zero sum of coordinates.\medskip

Let $K$ be the splitting field of the group algebra $\Q[G]$. Denote by $\X$
a full set of irreducible $K$-representations appearing in the decomposition of
$(K\Omega)^0$. Then the Galois group of the extension $K/\Q$ acts transitively on~$\X$.
This implies that the pair $(n_\chi,m_\chi)$ does not depend on $\chi\in\X$ where $n_\chi$
and $m_\chi$ are respectively the degree and multiplicity of~$\chi$. From this one can
deduce that the scheme of the group $G$ is pseudocyclic. Moreover, by the second part of
Theorem~\ref{pscyclic} this scheme is also commutative.

\subsection{The Terwilliger algebra of a pseudocyclic scheme}
Let $(\Omega,S)$ be an arbitrary scheme and let $T_\alpha(\Omega,S)$ be the Terwilliger
algebra of it at point $\alpha\in\Omega$, i.e. the subalgebra of $M_\Omega(\CC)$ generated
by $\CC S$ and the set of matrices $A(1_{\alpha s})$, $s\in S$. It immediately follows
that
$$
T_\alpha(\Omega,S)\subseteq\CC S_\alpha
$$
where $\CC S_\alpha$ is the adjacency algebra of the scheme $(\Omega,S_\alpha)$
(see Subsection~\ref{030909a}), and that equality holds exactly when the algebra in the
left-hand side is closed with respect to the Hadamard product.\medskip

Let $G=\aut(\Omega,S)$. As it was observed in \cite{IIY} it is always true that
$T_\alpha(\Omega,S)$ is always contained in the centralizer algebra $\CC G_\alpha$ of
the one-point stabilizer~$G_\alpha$. However, even for cyclotomic schemes this inclusion
can be strict. On the other hand, due to (\ref{210109c}) we have
$$
\CC S_\alpha\subseteq \CC G_\alpha.
$$
Thus from Theorem~\ref{220309h} and the above two inclusions it follows that
the Terwilliger algebra of the cyclotomic scheme of valency $k$ and rank at least
$4(k-1)k^3$ is coherent and coincides with $\CC G_\alpha$.

\subsection{Equivalenced schemes with bounded indistinguishing number}
Theorem ~\ref{220309h} may be strengthened if we replace the condition of being
pseudocyclic by bounding of an indistinguishing number of a scheme. Indeed,
the argument used in the proof of Theorem~\ref{220309h} yields us that a $k$-valenced
scheme with indistinguishing number $c$ of rank at least $ck^3$ is a separable Frobenius
scheme.

\section{Schemes, coherent configurations and permutation groups}\label{220309a}

\sbsnt{Definitions.}
Let $\Omega$ be a finite set and $S$ a partition of $\Omega\times\Omega$. Denote by $S^\cup$
the set of all unions of the elements of $S$. A pair $(\Omega,S)$ is called a {\it coherent
configuration} on $\Omega$ if the following conditions are satisfied:
\nmrt
\tm{S1} the diagonal $1_\Omega$ of $\Omega\times\Omega$ belongs to $S^\cup$,
\tm{S2} $S$ is closed with respect to $*$,
\tm{S3} given  $r,s,t\in S$, the number $c_{rs}^t=|\{\beta\in\Omega:\,(\alpha,\beta)\in r,\
(\beta,\gamma)\in s\}|$ does not depend on the choice of $(\alpha,\gamma)\in t$.
\enmrt
The elements of $\Omega$, $S$, $S^\cup$ and the numbers (S3) are called the {\it points},
the {\it basis relations}, the {\it relations} and the {\it intersection numbers}
of~$(\Omega,S)$, respectively. The numbers $|\Omega|$ and $|S|$ are called the {\it degree}
and {\it rank} of it. The unique basis relation containing a pair
$(\alpha,\beta)\in\Omega\times\Omega$ is denoted by $r(\alpha,\beta)$. The set
of basis relations contained in $r\cdot s$ with $r,s\in S^\cup$ is denoted by~$rs$.

\sbsnt{Homogeneity.}
The set $\Omega$ is the disjoint union of {\it fibers} of~$(\Omega,S)$, i.e. those
$\Delta\subset\Omega$ for which $1_\Delta\in S$. For any basis relation $s\in S$ there exist
uniquely determined fibers $\Delta,\Gamma$ such that $s\subset\Delta\times\Gamma$. Moreover,
the number $|\delta s|$ does not depend on $\delta\in\Delta$ and coincides with $c_{ss^*}^t$
where $t=1_\Delta$. We denote it by $n_s$. The coherent configuration $(\Omega,S)$ is called
{\it homogeneous} or a {\it scheme} if $\Omega$ is a fiber of it. In this case
$$
n_s=n_{s^*},\qquad s\in S,
$$
and the number $n_s$ is called the {\it valency} of $s$. We say that $(\Omega,S)$ is an
{\it equivalenced scheme of valency~$k$}, when $n_s=k$ for all $s\in S^\#$ where here
and below we put $S^\#=S\setminus\{1_\Omega\}$.

\sbsnt{Isomorphisms and schurity.}
Two coherent configurations are called {\it isomorphic} if there exists a bijection between
their point sets preserving the basis relations. Any such bijection is called an
{\it isomorphism} of these coherent configurations. The group of all isomorphisms of a coherent
configuration $(\Omega,S)$ contains a normal subgroup
$$
\aut(\Omega,S)=\{f\in\sym(\Omega):\ s^f=s,\ s\in S\}
$$
called the {\it automorphism group} of~$(\Omega,S)$. Conversely, let $G\le\sym(\Omega)$
be a permutation group and $S$ the set of orbits of the componentwise action of $G$
on~$\Omega\times\Omega$. Then $(\Omega,S)$ is a coherent configuration and we call it the
{\it coherent configuration of~$G$}. This coherent configuration is homogeneous if and only if
the group is transitive; in this case we say that $(\Omega,S)$ is the {\it scheme} of~$G$. A
coherent configuration on $\Omega$ is called {\it schurian} if it is the coherent configuration
of $2$-orbits of some permutation group on~$\Omega$.

\sbsnt{Algebraic isomorphisms and separability.}
Two coherent configurations $(\Omega,S)$ and $(\Omega',S')$ are called {\it algebraically
isomorphic} if
\qtnl{f041103p1}
c_{rs}^t=c_{r' s'}^{t'},\qquad r,s,t\in S,
\eqtn
for some bijection $\varphi:S\to S',\ r\mapsto r'$ called an {\it algebraic isomorphism}
from~$(\Omega,S)$ to~$(\Omega',S')$. Each isomorphism $f$ from~$(\Omega,S)$ to~$(\Omega',S')$
induces in a natural way an algebraic isomorphism between these schemes denoted by $\varphi_f$.
The set of all isomorphisms inducing the algebraic isomorphism~$\varphi$ is denoted by
$\iso(S,S',\varphi)$. In particular,
$$
\iso(S,S,\id_S)=\aut(\Omega,S)
$$
where $\id_S$ is the identical mapping on $S$. A coherent configurations $(\Omega,S)$ is
called {\it separable} if the set $\iso(S,S',\varphi)$ is non-empty for
{for each algebraic isomorphism~$\varphi$}.

\sbsnt{Semiregularity.}
A coherent configuration $(\Omega,S)$ is called {\it semiregular} if
\qtnl{040409b}
|\alpha s|\le 1,\qquad \alpha\in\Omega,\ s\in S.
\eqtn
A semiregular scheme is called {\it regular}; regular schemes are exactly thin schemes in the
sense of \cite{Zi1}. One can see that a coherent configuration (resp. scheme) is semiregular
(resp. regular) if and only if it is a coherent configuration (resp. scheme) of a semiregular
(resp. regular) permutation group. The proof of this statement as well as the next one can be
found in~\cite{EP00}.

\thrml{210109e}
Any semiregular configuration is schurian and separable.\bull
\ethrm

\sbsnt{One point extension.}\label{030909a}
Let $(\Omega,S)$ be a coherent configuration and $\alpha\in\Omega$. Denote by $S_\alpha$
the set of basis relations of the smallest coherent configuration on $\Omega$ such that
$1_{\alpha}\in S_\alpha$ and $S\subset S_\alpha^\cup$ (see~\cite{EP09}). The coherent
configuration $(\Omega,S_\alpha)$ is called the {\it $\alpha$-extension} (or a {\it one
point extension}) of $(\Omega,S)$. It is easily seen that
\qtnl{210109c}
\aut(\Omega,S)_\alpha=\aut(\Omega,S_\alpha).
\eqtn
Notice that the set $\alpha s$ is a union of some fibers of $(\Omega,S_\alpha)$ for all
$s\in S$, and the relation $t\cap (\alpha r\times \alpha s)$ belongs to the set
$S_\alpha^\cup$ for all $r,s,t\in S$.\medskip

Let $(\Omega',S')$ be a coherent configuration and $\varphi:(\Omega,S)\to(\Omega',S')$
an algebraic isomorphism. Let $\alpha'\in\Omega'$ and
$\psi:(\Omega,S^{}_{\alpha^{}})\to(\Omega,S'_{\alpha'})$ be an algebraic isomorphism
such that
\qtnl{290309z}
\psi(1_{\alpha^{}})=1_{\alpha'},\qquad \psi(s)\subset\varphi(\wt s),\ s\in S_\alpha,
\eqtn
where $\wt s$ is the unique basis relation of $(\Omega,S)$ that contains $s$. Then
$\psi$ is uniquely determined by $\varphi$. We say that $\varphi_{\alpha,\alpha'}:=\psi$
is the $(\alpha,\alpha')$-extension (or one point extension) of $\varphi$.

\sbsnt{$t$-condition.}\label{090209x}
The following definition goes back to \cite[p.70]{FKM}. Let
$(\Omega,S)$ be a coherent configuration. Two sets $\Delta,\Delta'\subset\Omega$ have
{\it the same type with respect to the pair $(\alpha,\beta)\in\Omega\times\Omega$} if
$\alpha,\beta\in\Delta\cap\Delta'$ and there exists a bijection
$\Delta\to\Delta',\delta\mapsto\delta'$ such that $\alpha'=\alpha$, $\beta'=\beta$ and
$$
r(\delta_1,\delta_2)=r(\delta_1',\delta_2'),\qquad \delta_1,\delta_2\in\Delta.
$$
Let $t\ge 2$ be an integer. The coherent configuration $(\Omega,S)$ satisfies the
{\it $t$-condition} at a relation $s\in S^\cup$ if for each $k=2,\ldots,t$ the number of
$k$-subsets of $\Omega$ of each fixed type with respect to the pair $(\alpha,\beta)\in s$ does
not depend on the choice of this pair. If the $t$-condition is satisfied at each $s\in S$, we
say that $(\Omega,S)$ satisfies the $t$-condition. It can be proved that the
scheme on $n$ points is schurian if and only if it satisfies the $t$-condition
for all $t=2,\ldots,n-1$.

\sbsnt{Indistinguishing number.}
Let $(\Omega,S)$ be a scheme. The {\it indistinguishing number} of a relation
$s\in S$ is defined to be the number
$$
c(s)=\sum_{t\in S}c_{tt^*}^s.
$$
The term goes back to \cite[p.563]{B81} where the number $n-c(s)$ with $n=|\Omega|$ was called
the {\it distinguishing number} of~$s$. Clearly, $c(1_\Omega)=n$. The maximum of $c(s)$,
$s\in S^\#$, is called the {\it indistinguishing number} of $(\Omega,S)$.

\lmml{180309u}
Let $(\Omega,S)$ be an equivalenced scheme of valency $k$. Then the arithmetical
mean of $c(s)$, $s\in S^\#$, equals $k-1$.
\elmm
\proof Set $n=|\Omega|$. Then $|S^\#|=(n-1)/k$ and $|s|=nk$ for all $s\in S^\#$.
Counting the number of $\Omega$-triples $(\alpha,\beta,\gamma)$
such that $r(\alpha,\beta)=r(\alpha,\gamma)\in S^\#$ by two ways we obtain that
$$
nk\sum_{s\in S^\#}c(s)=n(k-1)k|S^\#|
$$
whence the required statement follows.\bull

\sbsnt{Adjacency algebra.}\label{190809a}
Let $(\Omega,S)$ be a coherent configuration. The set $\{A(s):\ s\in S\}$ forms a
linear basis of an algebra $\CC S\subset\mat_\Omega(\CC)$. It is called the {\it adjacency algebra} of the coherent configuration $(\Omega,S)$.
From the definition it follows that it is closed with respect to the transpose and the
Hadamard product. In particular, it is semisimple. So by the Wedderburn theorem its
standard module $\CC\Omega$ is completely reducible. For an irreducible submodule $L$ of
$\CC\Omega$ corresponding to a central primitive idempotent~$P$ of the algebra~$\CC S$, we
set
\qtnl{190907}
n_\SP=\dim_\CC(L),\quad m_\SP=\rk(P)/n_\SP,
\eqtn
thus $m_\SP$ and $n_\SP$ are the {\it multiplicity} and the {\it degree} of the corresponding
irreducible representation of~$\CC S$. Clearly,
\qtnl{190309a}
|\Omega|=\sum_{P\in\P}m_\SP n_\SP,\qquad |S|=\sum_{P\in\P}n_\SP^2
\eqtn
where $\P$ is the set of central primitive idempotents of $\CC S$. The coherent
configuration $(\Omega,S)$ is called {\it commutative} if the algebra $\CC S$
is commutative (equivalently, $n_\SP=1$ for all~$P$).\medskip

Let $(\Omega,S)$ be a scheme.
It can be proved that $m_\SP\geq n_\SP$ for all~$P\in\P$ \cite{EP97}.
Besides, for the {\it principal} central primitive idempotent $P_0=\frac{1}{n}J$ where
$J$ is the all one matrix, we have $m_{\scriptscriptstyle{P_0}}=n_{\scriptscriptstyle{P_0}}=1$.
The corresponding one-dimensional irreducible representation of $\CC S$ takes $A$ to
$\tr(AA^*)$. In particular, we have
\qtnl{150309b}
n_rn_s=\sum_{t\in S}c_{rs}^tn_t.
\eqtn
We set $\P^\#=\P\setminus\{P_0\}$ and $\lg A,B\rg=\tr(AB^*)$ for all $A,B\in\CC S$.

\sbsnt{Intersection numbers.}\label{300109c}
There is a lot of useful identities for the intersection numbers of an arbitrary scheme \cite{Zi1}.
One of them is (\ref{150309b}), another one is
\qtnl{150309c}
n_tc_{rs}^{t^*}=n_rc_{st}^{r^*}=n_sc_{tr}^{s^*},\qquad r,s,t\in S.
\eqtn
We also need the following lemma.

\lmml{141208j}
Let $(\Omega,S)$ be a scheme and $r,s\in S^\#$. Then $c_{r^*s}^t\le 1$ for all
$t\in S$ if and only if $rr^*\cap ss^*=\{1_\Omega\}$.
\elmm
\proof We have
$$
n_rn_s\le n_rn_s+\sum_{t\in (rr^*\cap ss^*)^\#}c_{rr^*}^tc_{ss^*}^tn_t=
\lg rr^*,ss^*\rg
$$
and the bound is attained if and only if $rr^*\cap ss^*=\{1_\Omega\}$, where
$\lg rr^*,ss^*\rg=\lg A(r)A(r^*),A(s)A(s^*)\rg$. Similarly, from (\ref{150309b})
it follows that
$$
\lg r^*s,r^*s\rg=
\sum_{t\in S}(c_{r^*s}^t)^2n_t\ge
\sum_{t\in S} c_{r^*s}^tn_t=
n_rn_s
$$
and the bound is attained if and only if $c_{r^*s}^t\le 1$ for all $t\in S$. Since
$\lg rr^*,ss^*\rg=\lg r^*s,r^*s\rg$ we are done.\bull

\sbsnt{Frobenius groups.}
In this subsection we recall some well-known group theoretical facts on the
Frobenius groups~\cite{H59,Pa}. A non-regular transitive permutation group $G\le\sym(\Omega)$ is
called a {\it Frobenius group} if $G_{\alpha,\beta}=\{\id_{\Omega}\}$ for
all non-equal points $\alpha,\beta\in\Omega$. Any Frobenius group $G$ has a uniquely
determined regular normal subgroup $A$ called the {\it kernel} of $G$. Therefore
$G=AK$ where $K$ is a one point stabilizer of~$G$, and $\GCD(|A|,|K|)=1$.

\thrml{230309a}
Let $G\leq\sym(\Omega)$ be a finite {non-regular}
transitive permutation group with point stabilizer~$K$
and
$$
\theta=\sum_{\chi\in\Irr(G)} m_\chi \chi
$$
the decomposition of the permutation character $\theta$ of $G$ into irreducibles. Then
$G$ is a Frobenius group if and only if\, $\chi(1)/m_\chi=|K|$ for each
character $\chi\in\Irr(G)^\#$ with $m_\chi\ne 0$.
\ethrm
\proof To prove the necessity suppose that $G$ is a Frobenius group with a complement~$K$. Then
from \cite[p.318-319]{H59} it follows that given $\varphi\in\Irr(K)^\#$ the class
function $\chi_\varphi=\varphi^G-\varphi(1)\theta'$ is an irreducible character
of~$G$ of degree $\varphi(1)$ and
\qtnl{240309a}
\rho-|K|\theta'={\bf 1}+\sum_{\varphi\in\Irr(K)^\#}\varphi(1)\chi_\varphi
\eqtn
where ${\bf 1}$ and $\rho$ are the principal and regular characters of $G$, and
$\theta'=\theta-{\bf 1}$. Since $\chi_\varphi\in\Irr(G)$, we have
$[\rho,\chi_\varphi]=\chi_\varphi(1)=\varphi(1)$. Therefore by (\ref{240309a}) we
obtain
$$
\sum_{\chi\in\Irr(G)\setminus\Phi}\chi(1)\chi=|K|\theta'=\sum_{\chi\in\Irr(G)^\#}|K|m_\chi\chi
$$
where $\Phi=\{{\bf 1}\}\cup\{\chi_\varphi:\ \varphi\in\Irr(K)^\#\}$. Thus
$\chi(1)=|K|m_\chi$ for each character $\chi\in\Irr(G)^\#$ with $m_\chi\ne 0$.\medskip

To prove the sufficiency suppose that $\chi(1)/m_\chi=|K|$ for each non-principal
character $\chi\in\Irr(G)$ with $m_\chi\ne 0$. Since $\sum_\chi m_\chi^2=r-1$
where $r=|\orb(K)|$ (see \cite[Th.~16.6.14]{H59}), we have
$$
|\Omega|-1=\sum_{\chi\in\Irr(G)^\#}\chi(1)m_\chi=
\sum_{\chi\in\Irr(G)^\#}|K|m^2_\chi =
|K|(r - 1).
$$
Therefore $|K|=(|\Omega|-1)/(r-1)$. Since each non-trivial
orbit of $K$ has cardinality at most $|K|$ and there are $r-1$ orbits, we obtain that each
non-trivial $K$-orbit has cardinality $K$, that is $G_{\alpha,\beta}=\{\id_\Omega\}$
whenever $\alpha\neq\beta$.~\bull

\end{document}